\documentclass[leqno]{article}%
\usepackage{amsmath,amstext, amsthm, amssymb}
\newcommand{\modified}{}

\numberwithin{equation}{section}
\newtheorem{thm}{Theorem}[section]

\newtheorem{prop}[thm]{Proposition}
\theoremstyle{remark}
\newtheorem{rem}[thm]{Remark}
\theoremstyle{definition}
\newtheorem{definition}[thm]{Definition}
\newtheorem{example}{Example}[section]

\newcommand{\be}{\begin{equation}}
\newcommand{\ee}{\end{equation}}
\newcommand{\ba}{\begin{array}}
\newcommand{\ea}{\end{array}}

\renewcommand{\em}{\it}

\newcommand{\la}{\lambda} %
\newcommand{\eps}{\varepsilon}

\newcommand{\bea}{\begin{eqnarray}}
\newcommand{\eea}{\end{eqnarray}}

\newcommand{\Sum}{\sum_{n=0}^\infty}

\newcount\mins \newcount\hours \hours=\time \mins=\time
\divide\hours by60 \multiply\hours by 60 \advance\mins by-\hours
\divide\hours by60

\def\now{ \ifnum\hours>11 \ifnum\hours>12 \advance\hours by
-12 \fi \number\hours:\ifnum\mins<10 0\fi \number\mins\ pm,\ \else
\ifnum\hours=0 \hours=12 \fi \number\hours:\ifnum\mins<10 0\fi
\number\mins\ am,\ \fi}

\newcommand{\D}{{D}}

\title{%
Approximation
Operators, Exponential, $q$-Exponential, and Free Exponential Families}
    \author{W\l odzimierz Bryc\thanks{Research partially supported by NSF
grants \#INT-03-32062, \#DMS-05-04198} \\
Department of Mathematics \\
University of Cincinnati\\
Cincinnati, OH 45221-0025
\and Mourad E.H. Ismail%
\\ Department of Mathematics \\ University of Central Florida
\\  Orlando, FL 32816-1364
\\ }

\date{December 6, 2005}
\begin{document}

\maketitle

\tolerance=500

\begin{abstract}
Using the technique developed in approximation theory, we construct
examples of exponential families of infinitely divisible laws
which can be viewed as $\eps$-deformations of the  normal,
gamma, and Poisson exponential families.
Replacing the differential equation of approximation theory by a $q$-differential equation, we
define the $q$-exponential families, and we identify all $q$-exponential families with quadratic variance
functions when $|q|<1$. We  elaborate on the case of $q=0$ which is related to free convolution of measures.
We conclude by considering briefly the case $q>1$, and other related generalizations.
\end{abstract}

\noindent
\fbox{\begin{minipage}{4.5in}
{\it Running Title}.
Exponential Families \vskip2pt

\noindent
{\it Mathematics Subject Classification}. Primary: 62E10
Secondary:
41A35,
44A10,
33C45
\vskip2pt

\noindent {\it Key words and phrases}. Variance function,
exponential operators, $q$-Hermite polynomials, Al-Salam--Chihara
polynomials, free exponential families, freely infinitely divisible
laws, $q$-derivative, Hahn operator
\end{minipage}}
 \bigskip
\section{Introduction}
\subsection{Exponential Type Approximation Operators}
C. P. May %
\cite{May} introduced exponential type operators as
\bea S_\la(f)(m)
= \int_{\mathbb R}W_\la( m, u) \,f(u) du, \eea where $W_\la$ is a
generalized function satisfying the generalized differential
equation
\bea
\label{eqgendif} \frac{\partial W_\la}{\partial m} =\la W_\la
\frac{u-m}{v(m)}, \quad \la >0, \eea and $v$ is a polynomial of
degree at most 2. Moreover May assumed that $S_\la$ is a positive
operator and that \bea \label{normal} \int_{\mathbb R}W_\la( m, u)
\, du =1. \eea May knew the exact form of $S_\la$ in all possible
cases except when $v$ has two non-real zeros. May proved that
$S_\la$ and certain linear combinations of it approximate continuous
functions in the sense that $\lim_{\la \to \infty}S_\la(f, m) =
f(m)$. Later Ismail and May \cite{Ism:May} extended the
approximation theoretic study to the case when $v$ is
 an analytic and strictly positive function on
$(A,B)$, a component of $\{t:\; v(t)
>0\}$. They also identified $W_\la$ when $v(m)=1+m^2$, the general case with two complex roots.

In the above mentioned work, it was observed  that
\begin{equation}\label{Mean-variance}
\int_{\mathbb R}W_\la( m, u) \,u \,du =m, \qquad \int_{\mathbb R}
W_\la(m, u) \,(u-m)^2 du = \frac{v(m)}{\la}
\end{equation}
follow from \eqref{eqgendif} and \eqref{normal}. Hence $m$ and
$v(m)/\lambda$ are the mean and variance of $W_\la( m,u)$,
respectively.

The parameter $\la$ is important in approximation theory since as
$\la \to \infty$ the variance tends to zero and $W_\la$ becomes a
unit atomic measure concentrated at $u =m$.
Ismail and May \cite{Ism:May}  observed that the differential
equation \eqref{eqgendif} has at most one solution which satisfies
the normalization \eqref{normal} and makes $S_\la$ a positive
operator. They
 used the notation
\bea \label{eqdefq} q(m) = \int_c^m \frac{d \theta}{v(\theta)}, \;
\; c\in (A, B), \qquad g(q(m)) = q(g(m)) \equiv m. \eea Moreover
Ismail and May proved that \bea \label{formofS} S_\la(f,m) =
\int_{\mathbb{R}}C_\la( u) \exp\left(-\la \int_c^m
\frac{\theta-u}{v(\theta)}\, d\theta \right) f(u) \; du \eea and the
function (or generalized function) $C_\la(u)$ is computed by
inverting the Laplace transform \bea \exp\left(\la \int_c^m
\frac{\theta}{v(\theta)}\, d\theta \right)
=\int_{\mathbb{R}}C_\la(u) \exp\left( \la u \int_c^m
\frac{d\theta}{v(\theta)}\right) du. \nonumber \eea The above
formula is \bea \label{eqLaplace}
\begin{gathered}
\exp\left(\la \int_c^{g(z)} \frac{\theta}{v(\theta)}\, d\theta
\right) =\int_{\mathbb{R}}C_\la(u) \exp ( \la u z ) du,
\end{gathered}
\eea
and is valid for Re $z \in$ Range of  $q(m), m\in (A, B)$.
(Compare \cite[(2.1)]{Let:Mor}.) The theory of bilateral Laplace transform is in \cite{Zem}.

Ismail \cite{Ism} considered the case when $v(m)$ has a simple zero at an
end point,
which without loss of generality is taken as $m=0$.  %
He used the notation
\bea
\label{eqdefxieta}
\begin{gathered}
h(z) := \frac{1}{v(z)} - \frac{1}{z},\\
\xi = \xi(m) := \frac{m}{c} \exp\left\{\int_c^mh(\theta)d\theta
\right\}, \quad \eta(\xi) := m-c +\int_c^m \theta h(\theta) d\theta.
\end{gathered}
\eea He further assumed that $h(z)$  is analytic at $z=0$ and
$\eta'(0) \ne 0$. In his notation $W_\la$ is a discrete probability
distribution and takes the form \bea\label{eqdisoperform} W_\la(m,
du) = \Sum  \phi_n(\la) \, \exp\left(-\int_c^m\frac{\la\theta
-n}{v(\theta)}\, d\theta \right)\, \delta_{n/\la}(du), \eea where
$\{\phi_n:n =0, 1, \dots\}$ are generated by \bea\label{eqgenfun}
\exp(\la\eta(\xi)) = \Sum \phi_n(\la)\, \xi^n, \eea and
$\delta_a(du)$ is a unit atomic measure concentrated at
$a\in\mathbb{R}$.
Ismail also showed that $W_\la$ in \eqref{eqdisoperform} is
independent of the choice of $c\in(A,B)$.

\subsection{Exponential Families}
Fix a positive non-degenerate $\sigma$-finite measure $\mu$ on
$\mathbb{R}$ with the property that
$$
L(\theta)=\int_\mathbb{R} \exp(\theta u) \mu(du)<\infty
$$
for all $C<\theta<D$. Denote
$$\kappa(\theta)=\ln L(\theta).$$
The  exponential family generated by $\mu$ is the set of probability
measures
$$
\mathcal{F}(\mu):=\left\{P_\theta(du)=\exp(\theta u -
\kappa(\theta))\mu(du): \theta\in(C,D)\right\}.
$$
For a concise introduction, see \cite[Chapter 2]{Jorg}. Most authors
take $\theta$ from the largest admissible
 interval; ref. \cite{Let:Mor} restricts $\theta$ to a maximal open interval.

This family can be conveniently re-parameterized by the mean. Since
$\mu$ is non-degenerate, $\kappa(\cdot)$ is strictly convex
 so that $\kappa'()$ is strictly increasing on $(C,D)$;
it is also clear that $\kappa$ is analytic on $(C,D)$.
Let
\begin{equation}\label{A-B}
A=\lim_{\theta\to C^+} \kappa'(\theta),\; B=\lim_{\theta\to D^-}
\kappa'(\theta).
\end{equation}
Clearly, $\kappa':(C,D)\to(A,B)$ is invertible, and
$m=\kappa'(\theta)=\int_\mathbb{R} u P_\theta(du)\in(A,B)$. So  for
$\theta\in(C,D)$  probability measure $P_\theta$ is determined
uniquely by its mean $m\in(A,B)$. Let  $\psi$ be the inverse
function to $\kappa'$, i.e. $\kappa'(\psi(m))=m$ and
$\psi(\kappa'(\theta))=\theta$ for all $m\in(A,B)$,
$\theta\in(C,D)$. Then the  probability measures
\begin{equation}\label{P2W}
W(m,du):=P_{\psi(m)}(du), \; m\in(A,B)
\end{equation}
provide another parametrization of $\mathcal{F}(\mu)$.
Since $$
\int_\mathbb{R}uW(m,du)=m,$$
  this is parametrization by the means. The variance function
$V:(A,B)\to \mathbb{R}$ is now defined as
$$V(m)=\int (u-m)^2 W(m,du),$$
compare (\ref{Mean-variance}). Notice that $V(m)=\kappa''(\psi(m))$. It is known that  the variance
function $V$ together with $(A,B)$
 determines $\mu$ uniquely, see \cite[Theorem 2.11]{Jorg}, \cite[page 67]{Mor}, or \cite[Proposition 2.2]{Let:Mor}.

\subsection{Exponential Families and Exponential Operators}
The connection between exponential families and exponential operators has been noticed in \cite[Section 5]{diBuc:loeb},
see also \cite[Theorem 2]{Wedderburn74}.
Here we give a somewhat more precise version of this relation that allows for parameter $\la>0$ thus connecting
exponential operators with dispersion models \cite{Jorg}. \modified

 Suppose that a non-degenerate $\sigma$-finite measure $\mu$ with exponential moments of order $\theta \in (C,D)$
generates exponential family with the variance function $V(m)$,
$m\in(A,B)$. For natural $\la=1,2,\dots$ denote by $\mu_\la$ the
$\la$-dilation of the convolution power $\mu^{*\la}$, i.e.
$\mu_\la(U):=(\mu*\mu*\dots*\mu)(\la U)$. The natural exponential
family generated by $\mu_\la$   is the family of measures
$$
\mathcal{F}(\mu_\la):=\left\{P_{\la,\theta}(du)=\exp(-\theta u -
\kappa_\la(\theta)) \mu_\la(du):\; \theta\in(C\la,D\la)\right\},
$$
where $\kappa_\la(\theta)=\la \kappa(\theta/\la)$. In particular,
$\psi_\la(m)$ which is the inverse of $\kappa'_\la(\theta)$ is
$\psi_\la(m)=\la\psi( m)$ and the new variance function is
\begin{equation}\label{V-la}
V_\la(m)=\kappa_\la''(\psi_\la(m))=\frac{V(m)}{\la}.
\end{equation}
Notice that since  $\kappa'_\la(\theta)=\kappa'(\theta/\la)$, the
limits in \eqref{A-B} do not depend on $\la$. Parameterized
by the mean, the family is
$$\mathcal{F}(\mu_\la)=\{W_{\la}(m,du): \; m\in(A,B)\}.$$
We now verify that these measures satisfy equation (\ref{eqgendif}).
\begin{prop}\label{P1}
If a positive non-degenerate $\sigma$-finite measure $\mu$ with exponential moments of order $\theta\in(C,D)$ generates the natural
 exponential family with the variance function
$V(m)$ defined for $m\in(A,B)$, then for  natural $\la$ measures
$\mu_\la$ generate the natural exponential family $W_\la(m,du)$ such
that the corresponding integral operators
$$
S_\la(f)(m)=\int f(u)W_{\la}(m,du)
$$
are the exponential type operators which satisfy equation (\ref{eqgendif}) with  $\la=1,2,\dots$ and $v(m)=V(m)$ for $m\in(A,B)$.

\end{prop}
\begin{proof}
It is straightforward to verify that (\ref{eqgendif}) holds with $v(m)=V(m)$, $\la=1,2,\dots$.
Since
$$S_\la(f)(m)=\int f(u) \exp( \psi_\la(m) u - \kappa_\la(\psi_\la(m)))\mu_\la(du),$$
differentiating under the integral sign we get
\begin{multline*}
\int f(u)\frac{\partial}{\partial m}W_{\la,m}(du)\\ =\int f(u)
\psi_\la'(m)(u-\kappa_\la'(\psi_\la(m)))\exp\left( \psi(m) u -
\kappa(\psi(m))\right)\mu_\la(du).
\end{multline*}

As $\kappa_\la'(\psi_\la(m))=m$ and
$\psi_\la'(m)=1/\kappa''_\la(\psi_\la(m))=1/V_\la(m)=\la/V(m)$,
(\ref{eqgendif}) follows.
\end{proof}
\begin{rem}\label{R1}
If equation \eqref{eqgendif} has solution $S_\la(f,m)$ for
all $0<\la\leq 1$, and $m\in(A,B)$ then
 the exponential family generated by $\mu$  consists of infinitely divisible probability laws.
 \end{rem}
 \begin{proof}To prove infinite divisibility, without loss of generality we may concentrate on fixed $W_1(m_0,\mu)\in\mathcal{F}(\mu)$.
 It is well know that with the range of means $(A,B)$ kept fixed,  $\mathcal{F}(\mu)=\mathcal{F}(W_1(m_0,\mu))$,
 see \cite[Exercise 2.12]{Jorg}.

For $\la=1/k$ where $k=1,2,\dots$, let $W_\la(m,du), m\in(A,B)$  be the solution of (\ref{eqgendif}).
The variance function is $V(m)/\la=k V(m)$.
Denote by $\nu$ the dilation of measure $W_\la(m_0,du)$ by $k$.
By \eqref{V-la}, the exponential family $\mathcal{F}(\nu^{*k})$ has the same variance function $V(m)$ as the exponential family
$\mathcal{F}(W_1(m_0,du))$.
By uniqueness of   parametrization by the means, $W_1(m_0,du)=\nu^{*k}(du)$, so infinite divisibility follows.
 \end{proof}
Proposition \ref{P1} shows that  the celebrated result \cite[Section
4]{Mor} can be derived as a consequence of
 \cite[Theorem 3.3]{Ism:May}; the latter paper
contains also several cubic variance functions and other interesting
examples. Another interesting result \cite[Proposition 4.4]{Let:Mor}
is a consequence of \cite[Theorem 3.8]{Ism}.

\subsection{Notation}
We shall follow the terminology in \cite{Erd:Mag:Obe:Tri1}
for hypergeometric functions, namely that
\bea
\begin{gathered}
(a)_n :=1, \qquad (a)_n= \prod_{j=0}^{n-1}(a+j), \\
{}_2F_1\left( \left.
\ba{cc}
 a, b  \\
c
\ea
\right|z\right) := \Sum \frac{(a)_n \; (b)_n}{(c)_n\; n!}\; z^n.
\end{gathered}
\eea

The modified Bessel functions are \cite{Erd:Mag:Obe:Tri2}
\bea
\begin{gathered}
I_\nu(z) := \Sum \frac{(z/2)^{\nu+2n}}{n! \, \Gamma(n+\nu +1)},
\quad K_\nu(z) := \frac{\pi}{2} \,\frac{I_{-\nu}(z) - I_\nu(z)}{\sin
(\pi \nu)}.
\end{gathered}
\eea

The Lagrange expansion theorem \cite[(L), page 145]{Pol:Sze} says
that if $f(z)$, $\phi(z)$  are analytic in a neighborhood of $z=0$,
$\phi(0)\ne 0$ and $\xi := {m}/{\varphi(m)}$ then \modified \bea
\label{eqlagrange} f(m(\xi)) = f(0) + \sum_{n=1}^\infty \frac{\xi^n}{n!}
\, \left[\frac{d^{n-1}f'( x)[\phi(x)]^n}{dx^{n-1}}\right]_{x=0}.
\eea

By $1_{(a,b)}(u)$ we denote the indicator function of $(a,b)$.

 Occasionally, we also use the  $q$-notation
\begin{eqnarray*} (a;q)_n&:=&\prod_{k=0}^{n-1} (1-a
q^k),\\
(a;q)_\infty&:=&\prod_{k=0}^\infty (1-a q^k), \\
(a_1,a_2,\dots,a_m;q)_\infty&:=&(a_1;q)_\infty(a_2;q)_\infty\dots(a_m;q)_\infty,\\
{[n]_{q}} &:=&1+q+\dots +q^{n-1},\\ {[n]_{q}!} &:=&[1]_{q}[2]_{q}\dots
[n]_{q}=\frac{(q;q)_n}{(1-q)^n}, \\ %
\left[ \begin{array}{c}n \\ k
\end{array}\right] _{q}
&:=&\frac{[n]_{q}!}{[n-k]_{q}![k]_{q}!}=\frac{(q;q)_k(q;q)_{n-k}}{(q;q)_n},
\end{eqnarray*}
with the usual conventions $[0]_{q}=0,[0]_{q}!=1$.
Most of this notation is taken from  \cite{Gas:Rah}.
\section{Examples of Variance Functions}\modified
\subsection{$\eps$-Deformations of  Quadratic Variance Functions}
Letac and Mora \cite[page 3]{Let:Mor} raise the question of classifying
 exponential families with variances functions of the
form
\bea\label{Babel+}
P(m)+Q(m)\sqrt{R(m)},\eea
 where $P,Q,R$ are polynomials of degree
at most $3,2,2$ respectively. Letac \cite[page 74]{Let} \modified
initiated
the study of  variance functions \eqref{Babel+} when $P$ is a
multiple of $R$. The latter class was investigated by Kokonendji
\cite{kok95} who also gave an excellent overview of other known
cases.  Kokonendji  \cite{kok94}  used probabilistic techniques to
investigate variance functions in the Seshadri's class
$V(m)=\sqrt{R(m)}P(\sqrt{R(m)})$.
 This section  further advances the investigation of the variance functions  \eqref{Babel+}.

We use Proposition \ref{P1} to identify certain  exponential
families $\mathcal{F}_\eps$ with the variance
function of the form
\begin{equation}\label{V generic}
V(m)=(a m^2 + b m +c )\sqrt{1+\eps m^2}, \; \eps>0.
\end{equation}
These are $\epsilon$-deformations of the quadratic variance family
$\mathcal{F}_0$ analyzed in \cite{Ism:May} and \cite{Mor}. From
Mora's theorem \cite[Theorem 2.12]{Jorg}, as $\eps\to0$ while
$(A,B)$ is fixed, the corresponding probability laws in
$\mathcal{F}_\eps$  weakly converge to the respective
laws in $\mathcal{F}_0$.

We also give two examples of the  functions  \eqref{Babel+} which are not the variance functions.
\subsubsection{Continuous Exponential Families}\modified
  In this section we consider the following continuous $\eps$-deformations:
\begin{enumerate}
\item \label{case 1c}  $\eps$-Gaussian family $V(m)=(1+\eps m^2)\sqrt{1+\eps m^2}$,
\item\label{case 2c} $\eps$-gamma family $V(m)=m^2\sqrt{1+\eps m^2}$.
\end{enumerate}

The first case (\ref{case 1c}) gives an infinitely divisible family introduced in \cite{Barnd-78}, see
\cite[Example 2.5]{kok95} and
\cite[Exercise 3.2]{Jorg}.
\begin{thm}[{Kokonendji \cite{kok95}}]\label{Thm K} For $\la>0$, $\eps>0$, the exponential family with  the variance function
$$V(m)=\frac{1}{\lambda}(1+\eps m^2)\sqrt{1+\eps m^2},\; m\in\mathbb{R}$$
consists of the infinitely divisible probability laws with the
densities \bea \label{continii}
 \exp\left(\frac{\la}{\eps}\left(\frac{1+ u m \eps}{\sqrt{1+\eps m^2}}-1\right)\right)\,\frac{   \la }{\pi \eps \sqrt{1+\eps u^2}}\;
K_1( \frac{\la}{\eps}\sqrt{1+\eps u^2}) . \eea
\end{thm}

Before we give a  proof of Theorem \ref{Thm K} we show how we give a formal
argument. In the present case we have
$$
c=0, \;\; q(m) = \frac{m}{\sqrt{1+\eps m^2}}, \;\;  g(z) =
\frac{z}{\sqrt{1-\eps z^2}},
\;\;   %
\int_0^m \frac{\theta\, d\theta}{v(\theta)} = 1/\eps-
\frac{1}{\eps \sqrt{1+\eps m^2}}.
$$
Now \eqref{eqLaplace}, after $z \mapsto \sqrt{\eps} z/\la$, $u\mapsto u/\sqrt{\eps}$, and $\la\mapsto \la \eps$  becomes \bea
\label{eqltcaseii} \quad \exp\left(\la -  \sqrt{\la^2-z^2}\right)=
\int_{\mathbb R} \exp(u z) \, C_\la(u)du/\sqrt{\eps}, \quad {\rm Re}\; z \in
(-\la, \la). \eea If we know that the left-hand side of the above
equation is a bilateral Laplace transform we can use the inversion
theorem, Theorem 5a on page 241 of
Widder \cite[\S 6.5]{Wid}, %
 and see that
\bea
\begin{gathered}
e^{-\la} C_\la(-u)/\sqrt{\eps} = \frac{1}{2\pi i}\int_{-i\infty}^{i\infty}
\exp(-\sqrt{\la^2-v^2}) \, \exp(uv) \, dv \\
= \frac{1}{ \pi  }\int_{ -\infty}^{ \infty} \exp(-\sqrt{\la^2+v^2})
\, \cos(uv) \, dv.
\end{gathered}
\nonumber \eea Formula (26), page 16 of \cite{Erd:Mag:Obe:Tri4}
implies
\bea %
 C_\la (u) = \frac{\la\; e^{\la}\sqrt{\eps}}{ \pi \, \sqrt{1+u^2}}
K_1(\la\sqrt{1+u^2}). \eea

\begin{proof}[Proof of Theorem 2.2]
We verify  \eqref{eqltcaseii} directly. With the above $C_\la(u)$
the right-hand side of \eqref{eqltcaseii} is
\bea \label{2.5}%
\frac{\la\; e^{\la}}{ \pi}\int_0^\infty \cosh (uz)  \,
\frac{K_1(\la\sqrt{1+u^2})}{\sqrt{1+u^2}} \, du. \eea In view of
\cite[(7.2.40)]{Erd:Mag:Obe:Tri2}
$$
J_{-\frac12} (x) = \sqrt{\frac{2}{\pi x} } \; \cos x, \quad
K_{\frac12}(x) = \sqrt{\frac{\pi}{2 x} } e^{-x},
$$
we apply \cite[(7.14.46)]{Erd:Mag:Obe:Tri2} and conclude that the
expression in \eqref{2.5} equals the left-hand side of \eqref{eqltcaseii}.
Substituting back the original values of $\la,z,u$ we get \eqref{continii}.
\end{proof}

We now consider case (\ref{case 2c}), which yields the infinitely divisible distributions from
\cite[page 58]{Fel}.
\begin{thm}[{Letac \cite[page 46, Example 8.2]{Let}}]\modified
 For $\la>0$, $\eps>0$, the natural exponential family with the  variance function $$V(m)=\frac{m^2}{\lambda}\sqrt{1+\eps m^2}$$ defined on $m>0$,
 consists of the absolutely continuous infinitely divisible probability laws
\begin{multline}
\label{eps-gamma}
 \left(\frac{1+\sqrt{1+\eps m^2}}{\sqrt{\eps}m}\right)^{\la }\, \frac{  \la}{ u}\,
I_{ \la}(\sqrt{\eps} u) \exp\left(-\frac{\la  u (\sqrt{1+\eps m^2})}{m}\right)1_{(0,\infty)}(u)\,
du. \end{multline}
\end{thm}
\begin{proof}
We choose  $c=1/\sqrt{\eps}$ and apply \bea
\begin{gathered}
\int_c^m \frac{\theta}{v(\theta)} d\theta = \ln\left(\frac{m\,\sqrt{\eps} (1+\sqrt{2}\,)}{1+ \sqrt{1+\eps m^2}}\right), \\
\int_c^m \frac{d \theta}{v(\theta)} = \sqrt{2\eps} -
\frac{\sqrt{1+\eps m^2}}{m}.
\end{gathered}
\nonumber \eea Therefore \eqref{eqLaplace} gives
$$
\left(\frac{m\, \sqrt{\eps}(1+\sqrt{2}\,)}{1+ \sqrt{1+\eps m^2}}\right)^\la =
\int_0^\infty C_\la(u)\, e^{\la\, \sqrt{2\eps}\, u}
 \exp\left( -\la u (\sqrt{1+\eps m^2})/m\right)\, du.
$$
To invert the above Laplace we set $w =  (\sqrt{1+\eps m^2})/m$ so that
$m = 1/ \sqrt{w^2-\eps}$. Thus for $w >\sqrt{\eps}$ we need to invert
$$
\int_0^\infty C_\la(u)\, e^{\la\, \sqrt{2\eps}\, u}
 \exp( -\la u w)\, du =\eps^{\la/2}(1+\sqrt{2})^\la \, \left(w+ \sqrt{w^2-\eps}\right)^{-\la}.
$$
We use (28), page 240 in \cite{Erd:Mag:Obe:Tri4} to invert the above
Laplace transform and establish \eqref{eps-gamma}.
\end{proof}

\subsubsection{Discrete Exponential Families}\modified
In this section we consider the following cases:
\begin{enumerate}
\item\label{case 1d} the $\eps$-deformation  of the Poisson family $V(m)=m\sqrt{1+\eps m^2}$,
\item\label{case 4d}  the discrete $\eps$-deformation  of the Gaussian family $V(m)=\sqrt{1+\eps m^2}$.
\end{enumerate}

We first consider %
case (\ref{case 1d}).
 In this case $B = +\infty$ and we
choose $c=1/\sqrt{\eps}$. It is a calculus exercise to derive
$$ \int_1^m h(\theta) \, d\theta =
\ln \left(\frac{1+ \sqrt{2}}{1+\sqrt{1+m^2}}\right),$$
$$ \int_1^m
\theta h(\theta) \, d\theta = \ln \left(\frac{m+\sqrt{1+m^2}}{1+
\sqrt{2}}\right) +1-m.
$$
Hence \begin{multline}
 \xi(m) = \frac{(1+ \sqrt{2})\sqrt{\eps}m}{1+\sqrt{1+\eps m^2}}, \qquad
\eta(\xi(m)) =  \frac1{\sqrt{\eps}}\ln \left(\frac{\sqrt{\eps}m+\sqrt{1+\eps m^2}}{1+ \sqrt{2}}\right).
\end{multline}
With $\zeta(m) = \xi(m)/(1+\sqrt{2})$ it follows that $m =
\frac{2\zeta}{\sqrt{\eps}(1-\zeta^2)}$, so  that \bea \la \eta(\xi(m)) = \ln
\left(\frac{1+\zeta}{1-\zeta}\right)^{\la/\sqrt{\eps}}- \ln (1+ \sqrt{2})^{\la/\sqrt{\eps}}. \eea A
simple calculation shows that
$$\left(\frac{1+\zeta}{1-\zeta}\right)^\la = \Sum \frac{(\la)_n}{n!}
\; {}_2F_1\left( \left. \ba{cc}
-n, \quad \la\\
-\la-n+1 \ea \right|-1\right) \; \zeta^n.
$$
This proves the following theorem.
\begin{thm}[{Letac \cite[pg 98, (3)]{Let}}]\modified
 For $\la>0$, $\eps>0$, the exponential family
 with the  variance function $$V(m)=\frac{m}{\la}\sqrt{1+\eps m^2}$$
   defined on $m>0$,
 consists
of infinitely divisible discrete probability measures \bea \label{eqdiscretei}
\begin{gathered}
\left(\sqrt{\eps} m+\sqrt{1+\eps m^2}\right)^{-\la/\sqrt{\eps}} \Sum\frac{(\la/\sqrt{\eps})_n}{n!} \;
{}_2F_1\left( \left. \ba{cc}
-n, \quad \la/\sqrt{\eps}\\
-\la/\sqrt{\eps}-n+1 \ea
\right|-1\right) \\
\qquad \times  \left(\frac{\sqrt{\eps}m}{1+\sqrt{1+\eps m^2}}\right)^n
\delta_{n/\la}(du).
\end{gathered}
\eea
\end{thm}

We now consider %
{\bf Case (\ref{case 4d})}. This is again a known case:
\cite[Exercise 3.15]{Jorg} gives an answer in terms of the compound Poisson law,
\cite[Example 2.6]{kok95} writes the answer in terms of $\sum_{k\in Z} I_k(\la/\eps)\delta_k(du)$.
We remark that this is an example of a discrete indefinitely divisible
natural family to which  \cite[Proposition 4.4]{Let:Mor} or \cite[Theorem 3.3]{Ism} cannot be applied.

\begin{thm}[{Letac \cite[page 100, (8)]{Let}}]\modified
 For $\la>0$, $\eps>0$, the natural exponential family with the  variance function
 $$V(m)=\frac{1}{\lambda}\sqrt{1+\eps m^2}$$ defined on $m>0$,
consists
of infinitely divisible discrete probability measures
\begin{multline}
\label{eps-gauss}
e^{-\frac{\la}{\eps}\sqrt{1+\eps m^2}}\Sum \frac{\la^n}{(2\eps)^nn!} \sum_{k=0}^n \binom{n}{k}
\left(\sqrt{\eps} m+ \sqrt{1+\eps m^2}\right)^{ 2k-n}\delta_{(2k-n)\sqrt{\eps}/\la}(du).
\end{multline}
\end{thm}
\begin{proof}
We choose  $c=0$ and apply
\bea
\begin{gathered}
q(m)=\int_c^m \frac{d \theta}{v(\theta)} = \frac{\ln(\sqrt{\eps} m+\sqrt{1+\eps m^2})}{\sqrt{\eps}},\\
g(z)=q^{-1}(z)= \frac{\sinh (z \sqrt{\eps})}{\sqrt{\eps}},\\
 \int_c^m \frac{\theta}{v(\theta)} d\theta = \frac{\sqrt{1+\eps m^2}-1}{\eps}.
\end{gathered}
\nonumber
 \eea
 Therefore \eqref{eqLaplace} gives
$$ \int \exp(\la u z)C_\la(du)=\exp\left(\la(\cosh(\sqrt{\eps}z)-1)/\eps\right)=
\Sum e^{-\la/\eps} \frac{\la^n\cosh^n(\sqrt{\eps}z)}{\eps^n n!}.$$
Thus
$$
C_\la(du)=e^{-\la/\eps}\Sum \frac{\la^n}{(2\eps)^nn!} \sum_{k=0}^n \binom{n}{k}\delta_{(2k-n)\sqrt{\eps}/\la}(du)
$$
is just the
compound $\frac{1}{2}\left(\delta_{-\sqrt{\eps}/\la}+\delta_{\sqrt{\eps}/\la}\right)$-Poisson  law.
Using  the transform equation \eqref{formofS} we establish \eqref{eps-gauss}.
\end{proof}

\subsection{A Rational Variance Function}
Letac and Mora \cite[page 15]{Let:Mor} indicate that for $p_j>0$ the variance
function
$$V(m)=\frac{m}{(1-m/p_1)(1-m/p_2)\dots (1-m/p_k)}$$
corresponds to a discrete infinitely divisible exponential family which is difficult to determine explicitly.
Here we consider $v(m)=m/(1-m)$ which by dilation answers the question for $k=1$.

In this case  $\xi$ and $\eta$ of \eqref{eqdefxieta} %
with $c=1/2$ are
$$
\xi(m) = 2 \sqrt{e} me^{-m}, %
\quad \exp\left(\eta(\xi(m))\right) =
\exp\left(m-\frac{1}{2}\, m^2 - \frac{3}{8}\right)=\exp\left(-\frac{1}{2}\, (m-1)^2 + \frac{1}{8}\right).
$$
With $\phi(z) =  e^z/(2\sqrt{e}), f(m) = \exp\left(-\frac{\la}{2}(m-1)^2+\la/8\right)$ \modified
 in
\eqref{eqlagrange} we conclude that 
\begin{multline}
e^{\la \eta(\xi)} =e^{-3\la/8} +
 \sum_{n=1}^\infty  \frac{e^{\la/8}\xi^n}{2^n e^{n/2} n!}
 \left[\frac{d^{n-1}}{d x^{n-1}} e^{n x} \frac{d}{dx} \exp(-\la(x-1)^2/2)\right]_{x=0} \\
= e^{-3\la/8} +  e^{\la/8}\sum_{n=1}^\infty \frac{\xi^n}{2^n e^{n/2} n!} \sum_{k=0}^{n-1}
\binom{n-1}{k}\,n^{n-1-k}\left[\frac{d^{k+1}}{d x^{k+1}}
\exp\left(-\la(x-1)^2/2\right)\right]_{x=0}.
\end{multline}
 For $a>0$ we have
\bea \nonumber
\left[\frac{d^{k}}{d x^{k}}\, e^{-a(x-1)^2}\right]_{x=0}=e^{-a} a^{k/2} H_{k}(\sqrt{a}),
 \eea where $$H_n(x)=\sum_{j=0}^{\lfloor n/2\rfloor}\frac{(-1)^jn!}{j!(n-2j)!}\,x^{n-2j}$$ are Hermite polynomials.
Therefore
\eqref{eqgenfun} gives the following.

\begin{thm} For  $\la>0$, the natural exponential family with the  variance function
$$V(m)=\frac{m}{\lambda (1-m)}$$ defined on $0<m<1$,
 is generated by the  infinitely divisible discrete probability law $\mu_\la(du)=
 \Sum\phi_n(\la)\delta_n(du)$ with \modified
 \bea
\label{phicaseii}
\begin{gathered}
\phi_0(\la) :=\exp(-3\la/8), \qquad \qquad \qquad \\
\qquad \quad   \phi_n(\la) = \frac{e^{-3\la/8}}{2^n e^{n/2} n!} \sum_{k=0}^{n-1}
\binom{n-1}{k}\,n^{n-1-k}\left(\frac{\la}{2}\right)^{(k+1)/2}
H_{k+1}(\sqrt{\la/2}).
\end{gathered}\eea
\end{thm}
We note that  by \cite[Corollary
3.3]{Let:Mor} applied to the interval $M_F=(0,1)$ we have $\phi_n(\la) > 0$ for all $\la > 0$.

\begin{rem}
One can also write \eqref{phicaseii} as
$$
\phi_n(\la) = \frac{e^{-3\la/8}}{2^{n+1} e^{n/2} n!} \sum_{k=0}^{n}
\binom{n}{k}\,n^{n-k}\left(\frac{\la}{2}\right)^{(k-1)/2}
H_{k+1}(\sqrt{\la/2}).
$$
Similar calculations for $v(m)=m/(1+m)$ lead to $\phi_4(1)=-\sqrt{e}/64$, so this is not a variance function.
\end{rem}

\subsection{Positivity of $W_\la(m,du)$}
It is important to note that a given $v(m)$ does not
necessarily determine a distribution regardless of the choice of
$\la>0$.  Ismail gave such example in \cite{Ism}. In this section we
elaborate on this example and on another example  of the
form \eqref{Babel+}.

\begin{example}
Let $v(m)= m\sqrt{1-m}, m \in (0,1)$. With
$c = 1/2$ we find that \bea \xi(m) =
m\left[\frac{1-\sqrt{1-m}}{(\sqrt{2}-1)m}\right]^2, \quad \eta(\xi)
= \sqrt{2}-2\sqrt{1-m}. \nonumber \eea With $C := (\sqrt{2} -1)^2$
we have that
$$\eta(\xi) = \sqrt{2}- 2 +   \frac{4C\, \xi}{1+ C\, \xi}.
$$
Therefore \eqref{eqgenfun} becomes \bea \exp \left(\la (\sqrt{2}-2)+
4 \la C\xi/(1+C\xi) \right) = \Sum \phi_n(\la)\xi^n . \eea
 The information recorded so far is
from \cite{Ism}. Comparing (4.1) and (10.2.17), page 189 in
\cite{Erd:Mag:Obe:Tri2} we see that
\bea
\begin{gathered}
\phi_0(\la) = \exp\left(\la (\sqrt{2}-2)\right),\\
\phi_n(\la) = -4\la (-C)^n \exp\left(\la (\sqrt{2}-2)\right)\,
L_{n-1}^{(-1)} (4\la), \quad n >0,
\end{gathered}
\eea
where $L_n^{(-1)}(x)$ is the Laguerre polynomial.  Now
\eqref{eqdisoperform} shows that $W_\la$ is a probability
distribution if and only if $\phi_n(\la) \ge 0$ at the special value
of $\lambda$ under consideration. On the other hand Fej\'er's
formula \cite[Theorem 8.22.1]{Sze}
 shows that %
 $L_n^{(-1)}(4\la)$ is oscillatory at large $n$ for any
fixed positive $\la$. Thus there is no $\la$ for which $W_\la$ is a
probability distribution. This is an instance of the usefulness of
having the parameter $\la$.
\end{example}

\begin{example}
Let us now consider the case
$$
v(m) = \sqrt{1-m^2}.
$$
We take $c=0$. Thus $q(m) = \arcsin m, g(z) = \sin z$, and $\int_0^m
\frac{tdt}{v(t)} = 1- \sqrt{1-m^2}$. To determine $C_\la(u)$ we
need to invert
$$
\exp(\la(1-\cos z)) = \int_{\mathbb{R}} C_\la(u) \, e^{\la u z}\,
du,
$$
for all $z$, Re $z \in$ Range of $q(t)$,  $t\in (-\pi/2, \pi/2)$. Formula
(46) page 55 of \cite{Erd:Mag:Obe:Tri2} is \bea \label{eqKixform}
\int_0^\infty K_{ix}(a) \cos (xy) dx = \frac{\pi}{2}\, e^{-a\cosh
y}. \eea For large $p$ and fixed $a$, (19) page 88 in
\cite{Erd:Mag:Obe:Tri2} is
$$
K_{ip}(a) = \frac{\sqrt{2}\, \exp(-p\pi/2)}{(p^2-a^2)^{1/4}}(1+o(1)).
$$
Therefore \eqref{eqKixform} gives
\begin{equation}\label{K-wrong}
\frac{1}{\pi} \int_{\mathbb{R}}K_{ix}(\la) e^{ixy} dx  = e^{-\la \cos
y}.
\end{equation}
This implies \bea W_\la(m, u)= \frac{\la}{\pi} K_{i\la u}(\la)
\exp\left( \la \sqrt{1-m^2}  -\la u \arcsin m\right). \eea

Note that $ K_{i\la u}(\la)$ is real since $K_\nu(x) = K_{-\nu}(x)$ but it fails to be positive for any $\la>0$.
Indeed, the second derivative $\frac{d^2}{dy^2}$ of the right hand side of \eqref{K-wrong} at $y=0$ fails to be negative
as it equals $ \la e^{-\la}$.
\end{example}

\section{$q$-Exponential Families with $|q|<1$}
 Recall that for $-1<q<1$ the $q$-differentiation operator is
$$
(\D_{q,x}f)(x) := \frac{f(x)-f(qx)}{x-qx} \; \mbox{ for } x\ne 0.
$$

The  $q$-analogue of the differential equation \eqref{eqgendif} is
\bea \label{eqdefDq} \D_{q,m} w( m, u) = w( m, u)\,
\frac{u-m}{V(m)}. \eea
This equivalent to
\bea
\label{eqqt-V}
w(m,u) =  \frac{w(mq,u)}{1+m(1-q)(m- u)/V(m)}.
\eea
When $V(0) \ne 0$ we can rescale $m$ and $u$ by
a dilation to make $V(0) =1$. Now \eqref{eqqt-V} has the solution
\bea
\label{eqqW}
w(m,u) = C(u) \prod_{n=0}^\infty \frac{V(q^n m)}
{V(m) +m(1-q)(m- u)},
\eea
provided that the infinite products converge.

For compactly supported measures, the following extends the notion of exponential family from $q=1$ to $q\in(-1,1)$.
\begin{definition}\label{Def q-exp} A family of
probability measures
$$\mathcal{F}(V)=\left\{w(m,u)\mu(du): \; m\in(A,B)\right\}$$ is a $q$-exponential family with the variance function $V$
if \begin{enumerate}
\item $\mu$ is compactly supported,
\item $0\in(A,B)$ and $\lim_{t\to 0} w(t,u)=w(0,u)\equiv 1$ for all $u\in\mbox{\rm supp}(\mu)$,
\item $V>0$ on $(A,B)$ and \eqref{eqdefDq} holds for all $m\ne 0$.
\end{enumerate}
\end{definition}

Applying $\D_{q,m} $ to both sides of $\int w(m,u)\mu(du)=1$,  from \eqref{eqdefDq} we deduce
that
\begin{equation}\label{mean}
\int u w(m,u)\mu(du)=m.\end{equation}
This shows that family $\mathcal{F}(V)$ is parameterized by the mean.  Applying  $\D_{q,m}$ to both sides of \eqref{mean},
we deduce that
\begin{equation}\label{variance}
\int (u-m)^2 w(m,u)\mu(du)=V(m).
\end{equation}
Thus $V$ is  the variance function for $\mathcal{F}(V)$; compare \eqref{Mean-variance}.

We now show that quadratic variance functions determine $q$-exponential families uniquely.
\begin{thm}\label{T-q} If $\mathcal{F}(V)$ is a $q$-exponential family with the variance function $$V(m)=1+a m+b m^2$$
and $b>-1+\max\{q,0\}$
then
\begin{equation}\label{w-formula}
w(m,u)=\prod_{k=0}^\infty\frac{1+ a m q^k+ b m^2 q^{2k}}{1+(a-(1-q) u) m q^k+ (b+1-q) m^2 q^{2k}}
\end{equation}
and $\mu(du)$ is a uniquely determined probability measure with the absolutely continuous part
supported on the interval $\frac{a}{1-q}-\frac{2\sqrt{b+1-q}}{1-q}<u<\frac{a}{1-q}+\frac{2\sqrt{b+1-q}}{1-q}$
and no discrete part if $a^2<4b$
\end{thm}
We remark that for $b\geq 0$ the above family of laws $\mu$ appears in  \cite{Bryc-Wesolowski-03}
in connection to a quadratic regression problem. When $q\geq 0$,
one could also allow $b=-1/[N]_q$ for some integer $N\geq 1$ yielding a discrete
measure $\mu$ supported on $N+1$ points, compare \eqref{free mu} when $b=-1$.
\begin{proof}[Proof of Theorem \ref{T-q}] We rewrite
\eqref{eqdefDq} as \bea \label{eqq-normaleq} w( m, u) =
\frac{V(m)}{V(m)-(1-q)(u-m)m}w(qm, u). \eea
Thus
$$
 w( m, u) = w( q^{n+1} m, u)\prod_{k=0}^n \frac{1+ a m q^k+ b m^2 q^{2k}}{1+(a-(1-q) u) m q^k+ (b+1-q) m^2 q^{2k}}
$$
from which \eqref{w-formula} follows by taking the limit as $n\to\infty$.

We now recall that for $|t|$ small enough, \begin{equation}\label{w-sum}
w(t,u)=\Sum\frac{t^n}{[n]_q!}p_n(u),
\end{equation}
 is the generating function of the monic Al-Salam--Chihara polynomials
\begin{equation}\label{p-recur}
up_n(u)=p_{n+1}(u)+ a [n]_q p_n(u)+(1+b[n-1]_q)[n]_q p_{n-1}(u).
\end{equation}
This holds because the right hand side of \eqref{w-sum} satisfies \eqref{eqdefDq}, see \cite{Al-Salam76}.
Since $\mu$ is compactly supported,  we can integrate \eqref{w-sum} term by term for $|t|$ small enough;
we deduce that $\int p_n(u)\mu(du)=0$ for all $n\geq 1$. This determines probability measure $\mu$ as the
 measure of orthogonality of polynomials $p_n$.

Explicit formulas
 can be read out from
\cite[Chapter 3]{Askey-Ismail84MAMS}, see also \cite{Askey-Wilson85}. To use these results,
 we reparameterize \eqref{p-recur} as follows. Let $\widetilde{p}_n(x)=\alpha^{-n}p_n(\alpha x + \beta)$
 with $\alpha=\frac{\sqrt{b+1-q}}{\sqrt{1-q}}$, $\beta=a/(1-q)$. Then $\widetilde{p}_n(x)$ satisfy the three step
 recurrence
$$
(x -\widetilde{a} q^n)\widetilde{p}_n(x) =\widetilde{p}_{n+1}(x)+(1-\widetilde{b} q^{n-1}) [n]_q \widetilde{p}_{n-1}(x)
$$
with $\widetilde{a}=-\frac{a}{\sqrt{1-q}\sqrt{b+1-q}}$, $\widetilde{b}=\frac{b}{b+1-q}$.

\end{proof}
The
 technique we used in the proof of  will not work beyond polynomials of degree at most 2.
Al-Salam and Chihara \cite{Als-Chi87} proved
that the only orthogonal polynomials $\{p_n(x)\}$ with the generating function
\bea
\Sum p_n(x) t^n = A(t) \prod_{n=0}^\infty \frac{1-axH(tq^k)}{1-bxK(tq^k)}
\eea
where $A, H, K$ are formal power series with
\bea
A(t) = \Sum a_n t^n, \quad H(t) = \sum_{n=1}^\infty h_n t^n,
\quad K(t) = \sum_{n=1}^\infty k_n t^n
\eea
with $a_0h_1k_1 \ne 0$ and $|a| + |b| \ne 0$ are the
Al-Salam--Chihara polynomials if $ab =0$ and the $q$-Pollaczek
polynomials if $ab \ne 0$. Theorem \ref{T-q} corresponds to $a=0$.
The $q$-Pollaczek polynomials are in \cite{Cha:Ism87}.

(A related result appears also in \cite[Theorem 23]{Anshelevich03c}.)

\section{Free Exponential Families}
The case $q=0$ can be analyzed more directly. Since it is related to free
convolution of measures, it is of interest to elaborate explicitly on
the details.

\subsection{Free  Convolution of measures}
We  recall the analytic definition of the free  convolution  of compactly supported
probability measures due to Voiculescu \cite{Voiculescu86}, see also \cite[Section 2.4]{Voiculescu00}, \cite[Chapter 3]{Hiai-Petz00}.
The  Cauchy-Stieltjes transform \begin{equation}\label{Cauchy-Stieltjes}
G_\mu(z):=\int \frac{1}{z-u}\mu(du)
\end{equation}
 of a probability measure $\mu$ is analytic in $\Re z>0$. It is known that its inverse $G^{-1}(z)$ exists
 for $|z|$ large enough. The
  $R$-transform of $\mu$ defined as
$R_\mu(z)=G^{-1}(z)-1/z$ plays the role of the cumulant generating function.
 A probability measure $\mu$ is the free addictive convolution of probability measures $\mu_1,\mu_2$ if
$$R_\mu(z)=R_{\mu_1}(z)+R_{\mu_2}(z).$$
We write $\mu=\mu_1\boxplus \mu_2$.

The free cumulants of $\mu$ are the coefficients of  the expansions
\begin{equation}\label{free cumulants}
R_\mu(z)=\sum_{n=1}^\infty k_n(\mu)z^{n-1}.
\end{equation}

 \subsection{Exponential Families with  $q=0$}
As previously, we consider $A<0<B$ and assume that $V>0$ on $(A,B)$.

\begin{definition} A free exponential family with the variance
function $V(m)>0$ in a neighborhood of $0$ is a family of probability measures of
the form
\begin{equation}\label{W-free}
\mathcal{F}(V):=\left\{\frac{V(m)}{V(m)+m(m-u)}\mu(du):
m\in(A,B)\right\},
\end{equation}
where $\mu$ is a compactly supported probability measure.
\end{definition}
It is easy to verify that \eqref{W-free} defines a  family of measures which
fulfills all the requirements of Definition \ref{Def q-exp}, including  equation
\eqref{eqdefDq} with $q=0$.
It is also clear that the interval $(A,B)$ must be chosen so that the integral \eqref{W-free} converges.

For the purpose of determining measure $\mu$ alone, the role of the interval $(A,B)$ is insignificant. Namely, if
$V$ is a real analytic function at $0$, then
$\mu$ is determined uniquely by $V$. Indeed, since $V(0)\ne 0$,
the Cauchy-Stieltjes transform \eqref{Cauchy-Stieltjes}
is well defined for all real $z=m+\frac{V(m)}{m}$ large enough, i.e. for all $m$ close enough to $0$, and is given by
\begin{equation}\label{G2V}
G_\mu(z)=\frac{m}{V(m)}.
\end{equation}
This determines
$G_\mu(z)$ uniquely as an analytic function  outside of the support of $\mu$.

In particular, with  $V(m)=1+a m + b m^2$, equation $z=m+\frac{V(m)}{m}$ can be solved for
$m$, giving
$$
m=\frac{z-a - \sqrt{\left( a - z \right)^2-4\,\left( 1 + b \right)}}{2\,
         \left( 1 + b \right) } ,
             $$
             and
$$
G(z)= \frac{a + z + 2\,b\,z- {\sqrt{
         {\left( a - z \right) }^2-4\,\left( 1 + b \right)  }} }{2\,
     \left( 1 + a z\ + b\,z^2  \right) }.
$$
This Cauchy-Stieltjes transform appears in \cite[(2)]{Bozejko-Bryc-04}
in a non-commutative quadratic regression problem. It also appears
 in \cite[Theorem 4]{Anshelevich01}, \cite{Saitoh-Yoshida01}, and \cite[Theorem 4.3]{Bryc-Wesolowski-03}.
The corresponding laws are the free-Meixner laws
\begin{multline}\label{free mu}
\mu(du)=\frac{\sqrt{4(1+b)-(u-a)^2}}{2\pi (bu^2 +au +1)} 1_{(a-2\sqrt{1+b}, a+2\sqrt{1+b})} du
+p_1\delta_{u_1}+ p_2\delta_{u_2}.
\end{multline}
The discrete part of $\mu$ is absent except for the following   cases:
\begin{enumerate}
\item if $b=0, a^2>1$, then
$ p_1=1-1/a^2$, $u_1={-1/a}$, $p_2=0$.
\item if $b>0$ and $a^2>4 b$, then
$
p_1=\max\left\{0, 1-\frac{|a|-\sqrt{a^2-4b}}{2 b\sqrt{a^2-4b}}\right\}
$, $p_2=0$, and
$u_1=
\pm\frac{|a|-\sqrt{a^2-4b}}{2b}$ with the sign opposite to the sign  of $a$.
\item  if $-1\leq b<0$ then there are two atoms at %
$$
u_{1,2}=\frac{-a\pm \sqrt{a^2-4b}}{2b}, \; p_{1,2}= 1+\frac{\sqrt{a^2-4b}\mp a}{2 b\sqrt{a^2-4b}}.
$$
\end{enumerate}

This proves the free version of  \cite[Theorem 3.3]{Ism:May}, see also
\cite[Section
4]{Mor}.
\begin{thm}\label{T free quadr}
The free exponential family with the variance function $$V(m)=1+a m+b m^2$$
and $b>-1$
consists of probability measures \eqref{W-free} with $\mu$ given by \eqref{free mu}.

\end{thm}
We remark that if $b\geq 0$ then $\mu$ is infinitely divisible with
respect to  free convolution. In particular, up to dilation
and convolution with degenerate $\delta_m$ (i.e. up to "type")
measure $\mu$ is
\begin{enumerate}
\item  the Wigner's semicircle (free Gaussian) law if $a=b=0$; see \cite[Section 2.5]{Voiculescu00};
\item  the Marchenko-Pastur (free Poisson) type law if $b=0$ and $a\ne 0$; see \cite[Section 2.7]{Voiculescu00};
\item the free Pascal (negative binomial) type law if $b>0$ and $a^2>4b$; see \cite[Example 3.6]{Saitoh-Yoshida01};
\item the free Gamma type law if $b>0$ and $a^2=4b$;  see \cite[Proposition 3.6]{Bozejko-Bryc-04};
\item the free analog of hyperbolic type law if $b> 0$ and $a^2<4b$; see \cite[Theorem 4]{Anshelevich01};
\item the free binomial type law    if $-1\leq b <0$; see \cite[Example 3.4]{Saitoh-Yoshida01}.
\end{enumerate}

We conclude this section with the free version of Proposition \ref{P1}.
Recall that the
 $\la$-fold free convolution $\mu^{\boxplus \la}$
is well defined for the continuous range of values  $\la\geq 1$, see
 \cite{Nica-Speicher}. Let
$\mu_\la$ be the dilation by $\la\geq 1$ of the free convolution
power $\mu^{\boxplus \la}$.

\begin{prop}\label{P free}
If a compactly supported probability measure $\mu$  generates the free
 exponential family \eqref{W-free} with the real-analytic variance function
$V >0$ on $(A,B)$, then for  all $\la\geq 1$ measures
$\mu_\la$ generate the free exponential family with the variance function $V(m)/\la$.

Moreover, if for every $0<\la<1$ there is a $\mu_\la$ which generates the free exponential family
with the variance function $V(m)/\la$, then $\mu$ is infinitely divisible with respect to the free convolution.
\end{prop}
\begin{proof}
From \eqref{G2V} we determine
\begin{equation}\label{R}
R_\mu\left(\frac{m}{V(m)}\right)=m,
\end{equation}
 which determines the $R$-transform  $R_\mu$ uniquely in a neighborhood of $0$. Repeating the same calculation with the $R$-transform
 $R_{\mu_\la}$ of measure $\mu_\la$, we get
 $$R_{\mu_\la}\left(\frac{m}{V(m)}\right)=\la m.$$
 Thus
 $$
 \int \frac{V(m)}{V(m)+\la m(m-u)}\mu_\la(du)=1
 $$
 for all $|m|$ small enough, and $\mu_
 \la$ generates the corresponding free family \eqref{W-free}.

 The second part follows from the relation
 $$R_{\mu_\la}\left(\frac{m}{V(m)}\right)=\la R_\mu\left(\frac{m}{V(m)}\right),$$
which proves that $\mu$ is infinitely divisible with respect to the free convolution.
\end{proof}
\begin{rem} Combining \eqref{R} with \eqref{eqlagrange} we see that the  free exponential family with the
analytic variance function $V$ is defined by the unique centered probability measure $\mu$ with
 free cumulants
 $$k_{n+1}(\mu)=\left[\frac{1}{n!}\frac{d^{n-1}}{d t^{n-1}}V^n(t)\right]_{t=0},\; n=1,2,\dots.$$
\end{rem}

\section{$q$-Exponential Families for $q>1$}\label{Exp Fam q>1}
In this section it is convenient to set  $q =\frac{1}{p}$ with  $0<p<1$, and to use again auxiliary parameter $\la>0$.
The   $q$-analogue of the differential equation \eqref{eqgendif} is
\bea
\label{eqdefD1/q}
\D_{q,m} W_\la( m, du) =   W_\la(m, du)\, \frac{\la(u-m)}{v(m)}.
\eea
As previously, from $\int_{\mathbb{R}}W_\la(m, du)=1$ we deduce by $q$-differentiation that $\int_{\mathbb{R}} u W_\la(m, du)=m$
and
$\int_{\mathbb{R}}(u-m)^2 W_\la(m, du)=v(m)/\la$.

With
\bea \nonumber
\la_1:= \la(1-p)
\eea
we find that \eqref{eqdefD1/q} is
\bea
\label{eq1/q-normaleq}
W_\la(m, du) =  W_\la(p m, du)\left[1+ \frac{\la_1(u-pm)m}{v(p m)}\right].
\eea

We first consider \eqref{eq1/q-normaleq} when $q=\infty$. The general case of $v(0)\ne 0$ reduces to the case $v(0)=1$.
Substituting $p=0$ into \eqref{eq1/q-normaleq} we get
$$W_\la(m, du) =  \left(1+ \la  m u\right)C_\la(du),$$
which is  an analog of  equation \eqref{eqgendif} corresponding to  $q=\infty$.
From this, it is clear that
any probability measure
$C_\la(du)$ such that  $\int uC_\la(du)=0$, $\int u^2 C_\la(du)=1/\la$ determines its own $q$-exponential family
as long as $1+ \la u m \geq 0$ on the support of $C_\la(du)$.  Moreover, it is easy to see that the only
choice of $v(m)$ is a quadratic polynomial
$v(m)=1+\ b m-\la m^2$ where $b=la^2\int u^3 C_\la(du)$. This show that $q$-exponential families for $q=\infty$ are
not determined uniquely by their variance functions.

It is plausible that  non-uniqueness persists for all $q>1$. For example, the $q$-Hermite polynomials $\{h_n(x|q)\}$ in
\cite{Ism:Mas94}   correspond
to probability measures which are not determined uniquely by moments. The $N$-extremal solutions of the moment problem, \cite{Akh},
are given by a one-parameter family $\{\mu(du;a): a\in (p,1)\}$ which is completely characterized
by Ismail and Masson in
\cite{Ism:Mas94}, see also  Chapter 21 in \cite{Ism2005}. Unfortunately, the construction of the corresponding
 exponential family via equation \eqref{eqdefD1/q} led us to the family of measures
\bea
\label{eq1/q-normaldist}
W_\la(m, du) = \prod_{k=0}^\infty (1+\la_1 m u/q^k - \la_1 m^2/q^{2k} )
\mu(du;a)
\eea
 with negative densities.

The non-uniqueness within the class of quadratic variance functions $v(m)$
is confirmed by the following two examples.
\begin{example}\label{Example q-Laguerre 1}
Consider the absolutely continuous family with support in
   $(0,\infty)$ with the density
\bea
w_\la(m, u) = \frac{(p^{-\la}-1)^\la\, (p;p)_\infty \sin(\pi \la)}
{\pi m^\la \, (p^{1-\la};p)_\infty}\,
 \frac{u^{\la-1}}{(-u(p^{-\la}-1)/m;p)_\infty}.
\eea
This is the case of $p$-Laguerre polynomials \cite[\S 3.21]{Koe:Swa}.
With $q=1/p$, a calculation verifies that
\bea\label{D-Laguerre}
\D_{q, m} w_\la(m,u) = \frac{p(1-p^\la)}{m^2(1-p)}\,
w_\la(m,u) (u-m).
\eea
Now (3.21.2), page 108 of \cite{Koe:Swa} shows that
\bea
\int_0^\infty w_\la(m,u) \, du = 1,
\quad \int_0^\infty w_\la(m,u)\,  u \,  du = m.
\eea
(The latter integral follows also from the former  by $q$-differentiation and \eqref{D-Laguerre}.)
Thus $$\mathcal{F}=\{w_\la(m,u)1_{u>0}du: m>0\}$$ is parameterized by the mean.
Applying $\D_{q, m}$ again,
we get the variance function
\bea\label{*3}
V(m)=\int_0^\infty (u-m)^2\, W_\la(m,du)  = \frac{m^2}{\la_q}
\eea
with $\la_q=\frac{1-p}{p(1-p^\la)}$.
This is a %
continuous $q$-analogue of the gamma family with $v(m)=m^2$.
\end{example}

\begin{example}\label{Example q-Laguerre 2}
 For $m>0$, consider the family of discrete measures
 $$
 W_\la(m,du)=w_\la(m,u)\mu(du)
 $$
 with the density
 \bea
w_\la(m, u) =u^\la \frac{(-c,-p/c;p)}{(-cu,-cp^\la,-c^{-1}p^{-la+1};p)}, \; c=(p^{-\la}-1)/m
\eea
with respect to discrete measure
 $$\mu(du)=\frac{(p^\la;q)_\infty}{(p;p)_\infty}\sum_{n=-\infty}^\infty \delta_{p^n}(du).$$
This is again related to $p$-Laguerre polynomials \cite[\S 3.21]{Koe:Swa}.
With $q=1/p$, a calculation verifies that \eqref{D-Laguerre} holds.

Now (3.21.3), page 108 of \cite{Koe:Swa} shows that
\bea\label{*}
\int_{\mathbb{R}} w_\la(m,u) \, \mu(du) = 1.\eea
As previously, applying $\D_{q, m}$ to both sides of \eqref{*} and using \eqref{D-Laguerre} we get
\bea\label{*1}
\int_{\mathbb{R}} u w_\la(m,u) \, \mu(du) = m.\eea
Applying $\D_{q, m}$ to both sides of \eqref{*1} and using \eqref{D-Laguerre} again,
we get $V(m)= {m^2}/{\la_q}$, compare \eqref{*3}.
Thus  $\{W_\la(m,du):\; m>0\}$ is a
discrete $q$-analogue of the gamma exponential family;
it shares the variance function and the $q$-differential equation with the continuous  $q$-analogue of the gamma exponential family
 from the previous example.
\end{example}

\section{Shifted $q$-Exponential Families}
The special role played by $0$ in Definition \ref{Def q-exp} is due to the fact that $q$-derivative $\D_{q,x}$ is dilation
invariant but not translation invariant.
More generally,  we consider  the $L$-operator introduced by Hahn \cite{Hahn49}.
This is a $q$-differentiation operator centered at $\theta\in\mathbb{R}$, which we can write as
  \bea \label{Toperator}
({_\theta}\widetilde{\D}_{q, x}f)(x) = \frac{f(x) -
f(qx+(1-q)\theta)}{(1-q)(x-\theta)}, \; x\ne \theta, \; q\ne 1. \eea The usual $q$-derivative $\D_{q,x}$ corresponds to
$\theta=0$. For $\theta\ne 0$ a dilation reduces
all such operators  to  $\theta=1$, in which case we use shorter notation
\bea
\widetilde{\D}_{q, x}:={_1}\widetilde{\D}_{q, x}.
\eea
Note that
$$ \widetilde{\D}_{q, x}1=0,\; \widetilde{\D}_{q,
x}x= 1.$$
With $A\leq 1\leq B$ and $V>0$ on $(A,B)$, the shifted $q$-exponential family with  variance function $V$ is the family of probability measures
$$\mathcal{F}_\theta=\{w(m,u)\,\mu(du): m\in(A,B)\}$$
such that
$$
\widetilde{\D}_{q, m}w(m,u)=w(m,u)\frac{u-m}{V(m)}.
$$
In this discussion  we are less restrictive than in Definition \ref{Def q-exp}:
 in the admissible range of values of $\theta$ we include the end-points of $(A,B)$,
 and we allow non-compact  support for  $\mu$.
Such a generalization is beyond the scope of this paper, so we  give only one explicit example for
 $0<q<1$, $B=1$ and one for $q>1$, $A=1$.

\begin{example} Consider the case of the Wall polynomials, see
\cite[\S 3.20]{Koe:Swa}. In this case we have a family of discrete probability measures
$$
W(m,du)=w(m,u)\Sum \, \frac{q^n}{(q;q)_n}\,\delta_{q^n}(du),
$$
where the density is
\bea\label{3.20}
w(m,u) =  a^{\ln u/\ln q}(aq;q)_\infty
,\; a=(1-m)/q.
\eea
 From (3.20.3) in
\cite{Koe:Swa} we see that
$\int_{\mathbb{R}}p_1(u;a|q) w(m,u) du =0$, which implies
$\int_{\mathbb{R}}(1-ap  -u) w(m,u) du =0$. Hence the family
$$\mathcal{F}=\{W(m,du):\; 0<m<1\}$$ is again parameterized by the mean,
$$\int_{\mathbb{R}}u W(m,du)  = m.$$
From \cite[(3.20.3)]{Koe:Swa} we calculate the variance function
$$
V(m)=m(1-m)(1-q).
$$
Now \bea\label{1Dq bin} \widetilde{\D}_{q, m} w(m, u) =
\frac{u-m}{(1-q)m(1-m)}\, w(m,u). \eea Thus \eqref{3.20} defines a  shifted
analog of the $q$-Binomial family. Note that
although equation \eqref{1Dq bin} makes sense also for $q=0$,  it then gives a degenerated law $\delta_1$, not
 the translation of a free binomial law.
\end{example}

\begin{example} For $0<p<1$, let $q=1/p$ and consider the Al-Salam--Carlitz polynomials
$\{V_n^{(a)}(x;p)\}$, \cite[\S 3.25]{Koe:Swa}. Let
\bea\label{1-Poiss}
W(m,du)=w(m,u)\mu(du) = a^{-\ln u/\ln p} (aq/u;p)_\infty \mu(du),
\eea
where $\mu(du)=\Sum p^{n^2}/(p;p)_n \delta_{p^{-n}}(du)$
and $m = a+1$. Now with $q=1/p$ we have
$$
\widetilde{\D}_{q,m} w(m,u) = \frac{w(m,u)}{(m-1)(1-1/p)} \{1-
u(1-(m-1)/u)\},
$$
which simplifies to
\bea\label{V*}
\widetilde{\D}_{q,m} w(m,u) =  \frac{p}{1-p} \, w(m,u) \frac{u-m}{m-1}.
\eea
Since \cite[formula (3.25.2)]{Koe:Swa} implies $\int_{[0,\infty)} W(m,du)=1$,  therefore
applying $\widetilde{\D}_{q,m}$  and taking \eqref{V*} into account we deduce $\int_{[0,\infty)}
uW(m,du)=m$. Similarly $V(m) = (1-p)(m-1)/p$. Thus \eqref{1-Poiss} defines
 the family of measures
$$\mathcal{F}=\{W(m,du):\; m>1\},$$
which is  a
shifted $q$-analogue of the Poisson exponential family with $q=1/p>1$.
\end{example}

\subsection*{\bf Acknowledgements} The first named author (WB) thanks J. Wesolowski
for   suggesting the study of $q$-generalized exponential families, and additional references.
He also thanks   A. Vandal for a copy of his unpublished work and reference \cite{Wedderburn74}.
 and J. Chachulska for helpful comments. \modified

%
%


\end{document}